\documentclass[11pt]{article}

\usepackage[leqno]{amsmath}
\usepackage{amssymb}
\usepackage{amsthm}
\usepackage{amsfonts}
\usepackage{epsfig}
\usepackage{latexsym}
\usepackage{url}
\usepackage{ifthen}
\usepackage{comment}
\usepackage{subfigure}

\usepackage{natbib}
 \bibpunct[, ]{(}{)}{,}{a}{}{,}%

\newcommand{\field}[1]{\mathbb{#1}}
\def\R{\field{R}}                                

\def\Z{\field{Z}}

\def\uniond{\,{\buildrel d \over \bigcup}\,}

\def\eqd{\,{\buildrel d \over =}\,}

\def\blfootnote{\xdef\@thefnmark{}\@footnotetext}
\long\def\symbolfootnote[#1]#2{\begingroup%
\def\thefootnote{\fnsymbol{footnote}}\footnote[#1]{#2}\endgroup}

\newcommand{\Sscr}{{\cal S}}

\newcommand{\Pscr}{{\cal P}}

\newcommand{\Uscr}{{\cal U}}
\newcommand{\Bscr}{{\cal B}}
\newcommand{\Rscr}{{\cal R}}
\newcommand{\Mscr}{{\cal M}}

\newtheorem{theorem}{Theorem}

\newtheorem{corollary}{Corollary}
\newtheorem{lemma}{Lemma}
\newtheorem{assumption}{Assumption}
\newtheorem{property}{Property}

\newtheorem{example}{Example}

\newtheorem*{l3}{Lemma 3}

\newtheorem*{l5}{Lemma 5}

\DeclareMathOperator{\argmax}{argmax}

\newcommand{\eps}{{\epsilon}}
\newcommand{\Ascr}{{\mathcal A}}

\newcommand{\Hscr}{{\mathcal H}}

\renewcommand{\qed}{\hfill $\Box$ }

\usepackage{setspace}
\title{Stochastic Depletion Problems: Effective Myopic Policies for a class of Dynamic Optimization Problems}
\author{Carri W. Chan \footnote{Electrical Engineering, Stanford University, email:cwchan@stanford.edu}\\
\and
Vivek F. Farias \footnote{Sloan School of Management and Operations Research Center, Massachusetts Institute of Technology, email :vivekf@mit.edu}\\
}

\begin{document}
\maketitle
\def\thepage {} 
\thispagestyle{empty}

\begin{abstract}
This paper presents a general class of dynamic stochastic optimization problems we refer to as \emph{Stochastic Depletion Problems}. A number of challenging dynamic optimization problems of practical interest are stochastic depletion problems. Optimal solutions for such problems are difficult to obtain, both from a pragmatic computational perspective as also from a theoretical perspective. As such, simple heuristics are highly desirable. We isolate two simple properties that, if satisfied by a problem within this class, guarantee that a \emph{myopic} policy incurs a performance loss of at most 50 \% relative to the optimal adaptive control policy for that problem. We are able to verify that these two properties are satisfied for several interesting families of stochastic depletion problems and as a consequence identify efficient near-optimal control policies for a number of interesting dynamic stochastic optimization problems.

\end{abstract}

\newpage
\pagenumbering {arabic}

\section{Introduction}

This paper presents a general class of dynamic stochastic optimization problems we refer to as \emph{Stochastic Depletion Problems}. A number of challenging problems of practical interest are stochastic depletion problems. In fact, even special deterministic cases of some of these problems have been the focus of a great deal of recent research. Optimal solutions for such problems are difficult to obtain, both from a pragmatic computational perspective as also from a theoretical perspective. As such, simple heuristics are highly desirable. In this paper, we identify two simple properties that if satisfied by a stochastic depletion problem guarantee that a \emph{myopic} policy has expected value within a factor of $2$ of the optimal adaptive control policy for that problem. Such a myopic policy simply attempts to maximize one-period rewards (a notion we will formalize) and ignores information about the future evolution of the system making it practically attractive for many applications. We are able to verify that these two properties are satisfied for several interesting families of stochastic depletion problems and as a consequence identify efficient near-optimal control policies for a number of interesting dynamic stochastic optimization problems which we will highlight shortly.

Informally, a stochastic depletion problem is specified by \emph{item types} and \emph{activity sets}. The use of an activity results in the depletion of items of various types. In particular, the number of items of a particular type depleted at any point in time are randomly distributed according to a distribution specific to the activity employed at that time and the number of items of that type available. The sufficient statistics of these distributions are themselves specified by exogenous stochastic processes. New items of a given type may appear and existing items depart according to exogenous stochastic processes. Item depletion generates rewards, and activities must be selected adaptively over time to accomplish such depletion. An adaptive activity selection policy in this framework has knowledge of system dynamics and at any given time must select an activity. An optimal such policy generates maximum total expected reward.

We identify two simple properties we refer to as Value Function Monotonicity (VFM) and the Immediate Rewards (IR) property, that if satisfied by a stochastic depletion problem guarantee that a \emph{myopic} policy generates expected reward within a factor of $2$ of the optimal adaptive policy for that problem; that is, the myopic policy is a 2-\emph{approximation} algorithm for that problem. This policy is allowed to use all the information available up to the current point in time and maximizes expected reward earned over the following time step. Both properties are intuitive: the VFM property states that the optimal total expected reward (or value) accrued in the future starting from a particular state of the system is non-decreasing in the vector of available items at that state. In the other direction, the IR property states that the additional value gained by making available additional items at a particular state of the system is at most the reward earned for the depletion of those items.

We are able to verify both the VFM and IR properties for large families of stochastic depletion problems. These include stochastic depletion problems for which the total reward for items depleted over time is given by a non-decreasing submodular function of the vector of items depleted. An available item of a given type is depleted independent of all other items in the system with a probability that depends on time, the item type and the activity employed. We refer to these as Submodular Stochastic Depletion Problems. We are also able to address families of problems where the reward earned for depleting an item is a non-increasing function of the time of depletion, specific to that item's type. We refer to these as Decaying Linear Stochastic Depletion Problems. Our performance analysis of the myopic heuristic is sharp for both these families of stochastic depletion problems.

Our systematic study of the general class of problems presented in this paper results in a number of contributions which we now outline:

\noindent \textbf{Stochastic Control Problems: } Our framework lets us easily recognize and analyze simple near-optimal policies for several high-dimensional stochastic control problems for which finding optimal control policies is otherwise difficult.

For instance, control problems pertaining to several interesting discrete-time queueing models with general arrival processes and geometric service times may easily be reduced to stochastic depletion problems provided one allows for service disciplines with pre-emption. As an example, we consider a discrete time equivalent of a well studied `call-center' queueing model that has been an important subject of recent research (see for instance \cite{Harrison, Bassamboo_QSTA, Bassamboo_OR}). Finding optimal service policies for such problems is typically challenging. We identify a simple myopic policy that is near-optimal for a broad interesting class of performance metrics. Our policy is analogous to the so-called `$c-\mu$' scheduling rules (see for instance, \cite{VanMieghem}).

As another example we consider problems of dynamic `product line design'. These are problems where a firm must dynamically adjust the assortment of products it offers for sale so as to maximize expected revenues. Sales of a given product are influenced by the entire assortment of products offered as well as the prevailing market sizes for various customer segments. This represents an important generalization of static product line design problems of the type considered by \cite{VanRyzinMahajan99}. Our analysis yields the surprising conclusion that the optimal policy is well approximated by solving a sequence of \emph{static} product line design problems.

\noindent\textbf{Online algorithms for stochastic variants of well studied deterministic problems: }
We show that a myopic policy earns expected rewards that are within a small constant factor of the optimal adaptive policy for what we believe to be important online stochastic generalizations of a number of problems studied typically in deterministic settings. The stochastic generalizations we present go beyond what may be modeled in the traditional online versions of these problems and incorporate features we view as highly desirable from a modeling perspective. We present approximation guarantees for these generalizations that are typically no worse (and sometimes, better) than the best known guarantees for their deterministic counterparts.

For instance, we are able to provide an efficient myopic policy that is a 2-approximation algorithm for a stochastic broadcast scheduling problem. Successful data transmission in our broadcast scheduling model is \emph{stochastic} which makes it naturally applicable in several communications engineering contexts. With a bound on the maximal number of simultaneous broadcasts the best known scheme for \emph{the deterministic} broadcast scheduling problem (where page transmission is successful with probability $1$) is an offline $4$-approximation (\cite{BarNoy02}). With no constraints on the maximal number of simultaneous broadcasts, the best known offline algorithm for the deterministic broadcast scheduling problem is a $4/3$-approximation due to \cite{Gandhi02} and the best known online algorithm for the same scenario is a $2$-approximation due to \cite{Kima04}. Our contribution to this body of work is to demonstrate that even with stochasticity in data transmission one can still achieve good constant factor performance guarantees via an adaptive myopic algorithm.

As another example, we consider a stochastic generalization to the AdWords Assignment problem (\cite{Fleischer}) where revenues are generated via stochastic clicks on placed ads as opposed to simply via the placement of an ad as in the deterministic version of the problem. Using an appropriate approximation algorithm as a subroutine for the myopic problem we demonstrate an efficient $3+\eps$ myopic scheme which matches the best known $3+\eps$-approximation algorithm for the deterministic version of the AdWords Assignment problem due to \cite{Goundan}.

\noindent \textbf{A unified framework: } We provide an elegant, unified framework for the design and analysis  of stochastic optimization problems that is analogous to that for submodular maximization over simple matroids. In particular, submodular maximization problems over several simple matroids (such as the cardinality and partition matroid) can be captured as submodular stochastic depletion problems. A number of interesting problems (such as the AdWords Assignment problem) are known to be examples of such submodular optimization problems (see \cite{Goundan} for example). The stochastic depletion framework provides a natural vehicle for the analysis of stochastic variants of such problems wherein the notion of selecting a set element translates to an attempt at selection; the success of such an attempt is specified by an exogenous stochastic process.


We believe that the characterization of dynamic stochastic optimization problems that admit simple near-optimal control policies provided in the present work is likely to allow for a simple analysis of many problems beyond the handful of examples we have alluded to above. In particular, the two abstract properties that guarantee the effectiveness of a myopic policy are typically not hard to recognize and could potentially be established for families of problems outside those discussed here. The remainder of this paper is organized as follows: In Section \ref{sec:model} we formally specify the class of stochastic depletion problems. Section \ref{sec:proofs} presents a myopic heuristic for stochastic depletion problems and identifies two simple properties - the VFM and IR properties - that if satisfied by a stochastic depletion problem guarantee that a myopic policy is a 2-approximation algorithm for that problem. Section \ref{sec:families} verifies the VFM and IR properties for two general families of stochastic depletion problems - Submodular Stochastic Depletion problems and Linear Decaying Stochastic Depletion problems while the following two sections discuss a number of applications that lie within these families. Section \ref{sec:conclusion} concludes with a perspective on interesting directions for further work.

\section{Model}
\label{sec:model}
We are given a collection of \emph{items}, each of which belongs to one of
$M$ types indexed by $m$. There can be at most $\bar{x}_m$ items of
type $m$ available at any time. Items are
depleted via the execution of a suitable activity from a set of
feasible activities $\Ascr$, and depletion of a set of items garners a non-negative reward we will formalize shortly. Time is discrete (indexed by $t \in [0,T]$) and in each time step one must
choose to employ some activity from $\Ascr$; we let $T$ denote the length of the time horizon. We let $i$ index the
elements of $\Ascr$ and denote a general element of $\Ascr$ by $A$. Let $x_{t,m}$ denote the number of items of type
$m$ that remain at the start of the $t$th time-step. Assuming one
chooses action $A \in \Ascr$ in the $t$th time-step, the number of items of each type $m$ depleted within that time-step is given by a $\times_{m=1}^M \{0,\dots, x_{t,m}\}$-valued random vector, $X_{t}^A$, whose \emph{sufficient statistics} at time $t$ are specified by an exogenous stochastic process $\{P_{t}(A)\}$ specified for all $A$ and taking values in some compact set $\Rscr$. $X_{t}^A$ is assumed independent of the past given $P_{t}(A)$ and the vector of non-depleted items $x_{t}$. For example, we may have that $\{P_{t}(A)\}$ is a $[0,1]^M$-valued stochastic process and assuming one chooses action $A$ at time $t$, $X^A_{t,m}$ is an independent Binomial-$(x_{t,m},P_{t,m}(A))$ random variable. In what follows we drop the superscript $A$ from $X_t^A$ for economy, as the dependence on $A$ will be clear from context.
We have $x_{t+1,m} = x_{t,m} - X_{t,m}$ for all $m$, and receive a total
reward of $g(x_t,x_{t+1},t)$ where $g:\Z_+^M \times \Z_+^M \times [0,T]\rightarrow \R_+$ satisfies:
\begin{assumption}
\label{as:g}
For all $x,x' \in \Z_+^M$, $g(x,x',\cdot)$ is a non-increasing, non-negative function. In addition, we assume $g(x,x',T) = 0$ for all $x,x' \in \Z_+^M$.
\end{assumption}

Our objective is to design an adaptive scheduling
policy that maximizes total expected reward earned within the first
$T$ time-steps.
We define as our \emph{state-space} the set

\[
\begin{split}
\Sscr = &  \Big\{ (x,t,p_1,p_2,\dots,p_{|\Ascr|}):
x \in \times_m \{0,1,\dots,\bar{x}_{m}\}, 0 \leq t \leq T, p_i \in \Rscr^{t+1} \ \forall i \Big\}.
\end{split}
\]

In particular, a state is associated with a vector of items remaining to be depleted, time and a history of the $P_t$ processes. We denote by $x(s)$ the projection of $s$ onto its first co-ordinate and similarly employ the notation $t(s)$, and $p_i(s)$ for $i=1,2,\dots,|\Ascr|$. We let the random variable
$S_t \in \Sscr$ denote state in the $t$th epoch.

Finally, we define the random reward function $R: \Sscr \times \Ascr
\rightarrow \R_+$ according to $R(s,A) = g(x(s), x(s) - X_{t(s)},t(s))$,
where $X_{t(s)}$ is a random vector with sufficient statistic $P_{t(s)}(A)$.
Now since at time $t$, the realization of $X_t$ from taking a particular action is unknown, \emph{any} control policy is a-priori unaware of the exact reward accrued from a particular action. Only the statistics of this reward are known. More specifically, a control policy $\pi$ is a mapping from $\Sscr$ to the set of feasible activities $\Ascr$, and we denote by $\Pi$ the set of all
such policies.  Define the expected total reward-to-go under a
policy $\pi$ starting at state $s$ according to:
$$J^\pi(s) = E\left[\sum_{t'=t(s)}^{T-1} R(S_{t'},\pi(S_{t'})) | S_{t(s)} = s \right].$$
We let $J^*(s) = \max_{\pi \in \Pi} J^{\pi}(s)$ denote the expected
total reward-to-go under the optimal policy $\pi^*$ given by $\pi^*(s) \in \argmax_{\pi
\in \Pi} J^{\pi}(s)$. We will refer to the problem of finding such
an optimal policy $\pi^*$ as a \emph{Stochastic Depletion Problem}.

We remark that our formulation permits modeling exogenous item arrivals and deadlines on the latest time of depletion for a given item. In particular, assuming without loss that $\bar{x}_m = 1$, that is, a given item type can have at most a single item (otherwise, we could simply refine the definition of a type), we associate with each type an arrival time $\tau_m$ and a deadline $d_m \geq \tau_m$. One may then assume that the $P_t$ processes are such that $X_{t,m} = 0{\rm \ a.s.}$ for all $A \in \Ascr$ if $t \notin [\tau_m,d_m]$ in order to model the fact that item type $m$ arrives at time $\tau_m$ and may not be depleted beyond time $d_m$; see Section \ref{ssec:broadcast} for a concrete illustrative example. Such a formulation succinctly assumes a (known) bound on the total number of arrivals in any given period.

The optimal reward-to-go function (or \emph{value} function) $J^*$ and the optimal scheduling
policy $\pi^*$ can in principle be computed via dynamic programming: In particular, letting $S(s,A)$ denote the random next state encountered upon employing activity $A$ in state $s$ define the
dynamic programming operator $\Hscr$ according to:
\begin{equation}
\label{eq:DPoperator} (\Hscr J)(s) = \max_{A \in \Ascr}
E\left[R(s,A) + J(S(s,A)) \right].
\end{equation}
for all $s \in \Sscr$ with $t(s) \leq T-1$. $J^*$ may then be found as the solution to the Bellman equation
$\Hscr J = J$, with the boundary condition $J(s') = 0$ for all $s'$
with $t(s') = T$. The optimal policy $\pi^*$ may be found as the
greedy maximizer with respect to $J^*$ in \eqref{eq:DPoperator}.
Of course, this approach is computationally intractable: even in the event that the $P_t$ processes are known a-priori, the state space (the set of all $(x,t)$) is exponentially large. As such, this makes solution of a general stochastic depletion problem pragmatically difficult.

In addition to the above informal description for why we might
expect finding an optimal solution to be a difficult task, one may easily see that special classes of stochastic depletion problems are NP-hard. We consider one such class here for completeness: Consider `rational clairvoyant' stochastic depletion problems where the $P_{t}$
sequences are rational valued deterministic sequences in $[0,1]^M$. $X^A_{t,m}$ is a Binomial-$(x_{t,m},P_{t,m}(A))$ random variable for all $t,m,A$, independent of the past and $X_{t,m'}^{A'}$ for all $m' \neq m, A' \neq A$. We assume $g(x_t,x_{t+1},t) = \sum_m w_m (x_{t,m} - x_{t+1,m})$ where $w_m \geq 0$ for all $m$. The input to such a problem consists of $|\Ascr|$ rational-valued
sequences of length $T$, and $2M$ rational numbers representing the
initial number of jobs of each type and the reward constants, $w_m$,
for each type. One may then construct a polynomial time reduction
from the set-cover decision problem (which is NP-complete) to the
rational clairvoyant stochastic depletion problem. This is formalized in the
following result:

\begin{theorem}
The rational clairvoyant Stochastic Depletion problem is NP-hard.
\end{theorem}
\begin{proof}
We reduce the Set-Cover decision problem to the rational clairvoyant Stochastic
Depletion problem. An instance of the Set-Cover decision problem is specified
by a ground set $\Uscr$, cover set $\Bscr \in 2^\Uscr$ and an integer $k$ (where $k \leq |\Uscr|$ without
loss) and we must decide whether a cover (that is a subset of
$\Bscr$ whose union is a superset of $\Uscr$) of size $\leq k$
exists. We reduce this question to the optimal solution of the
following rational clairvoyant stochastic depletion instance: We consider a problem with
$|\Uscr|$ item-types, and assume we have a single item of each type.
We let each set in $\Bscr$ correspond to a feasible activity in the
sense that the use of that activity results in the depletion of all
items in that set with probability $1$ and the depletion of
item outside that set with probability $0$ in any time slot.
We let the depletion of a single item result in unit reward, $w_m =
1$, and assume that the time horizon for scheduling is $k$. Assuming
a polynomial time algorithm for rational clairvoyant stochastic depletion, the
reduced problem would require time that is ${\rm poly}(|\Uscr|,
|\Bscr|,k) = O({\rm poly}(|\Uscr|, |\Bscr|))$. If the optimal
solution to this instance of the stochastic depletion problem has
total reward $|\Uscr|$ we know that there exists a set cover of size
$\leq k$. Conversely, if there exists a set-cover of size $\leq k$,
then there exists a depletion policy with total reward $|\Uscr|$. Our
reduction is thus many-one and polynomial in the size of the input.
This completes the proof.
\qed
\end{proof}

We remark that the above reduction can also be used to reduce an interesting \emph{optimization} problem related to set covering - namely that of maximum set coverage where one may pick at most $k$ elements of $\Bscr$ so as to cover as many elements of $\Uscr$ as possible - to a stochastic depletion instance. In fact, this reduction is a special class of a useful set of reductions we will explore in Section \ref{sec:matroid}.

We have, in this section, introduced a general class of dynamic stochastic optimization problems that as we shall see in later sections admit a number of interesting applications. Computing optimal solutions for such problems is evidently hard; the next section will
present and analyze a natural, simple to implement heuristic for such problems.

\section{A Myopic
Heuristic for Stochastic Depletion}\label{sec:proofs}

A natural heuristic policy one may consider for a stochastic depletion problem is given by the myopic policy which in state $s$ chooses an activity set $A$ that
maximizes expected reward earned over the following time-step. That is
$$
\pi^{g}(s) \in \argmax_{A \in \Ascr} E[R(s,A)].
$$
Such a policy is adaptive but ignores the evolution of the system and the impact of the present choice of activity on rewards in future states. The set $\Ascr$ in the myopic problem above is potentially exponentially large. In many cases however, this set has an implicit polynomial sized representation (for instance, $\Ascr$ may correspond to a set of matchings) and the myopic maximization problem is efficiently solved. We will later also address the case where the myopic maximization problem is difficult but one has access to an appropriate near-optimal oracle.

\subsection{The Myopic Heuristic is an online $2$-Approximation Algorithm: Sufficient Conditions}\label{ssec:bound}

Our objective in this section will be to identify stochastic depletion problems for which the myopic heuristic is near-optimal. In particular, we will identify stochastic depletion problems for which we will have for any state $s \in \Sscr$, $J^*(s)/J^{\pi^g}(s) \leq 2$.

Noting that the myopic heuristic does not utilize any information about the evolution of the $P_t$ processes, we will simply assume that these are a-priori given sequences. In particular, we will compare the performance of the myopic heuristic to that of an optimal \emph{clairvoyant} algorithm that knows the realizations of the $P_t$ processes a-priori. Since an optimal clairvoyant policy must dominate the optimal policy, it will suffice to demonstrate performance guarantees relative to the optimal clairvoyant policy. Such an optimal clairvoyant policy may be computed over a reduced state-space:
\[
\begin{split}
\Sscr =  \Big\{ (x,t): x \in \times_m \{0,1,\dots,\bar{x}_{m}\}, 0
\leq t \leq T\Big\}.
\end{split}
\]
In the sequel, we will only consider such clairvoyant optimal policies; any reference to an optimal policy or value function in the sequel will pertain to an optimal policy or value function for the clairvoyant problem. Comparing performance to a clairvoyant policy yields performance guarantees that are valid over individual sample paths of the $P_t$ processes. In particular, our guarantees will imply that the myopic heuristic is a \emph{2-competitive online algorithm} where the optimal scheme is allowed knowledge of entire sample-paths of the $P_t$ processes but \emph{does not know} the realization of $X_t$ until time $t+1$ (which is somewhat different from the typical competitive analysis setting).

We now identify two properties that if satisfied by the optimal clairvoyant value function $J^*$, will imply our desired approximation guarantee.

\begin{property}
\label{prop:vfm}
\noindent \textbf{Value Function Monotonicity: } Consider states $s,s'$ satisfying $x(s) \geq x(s'), t(s) = t(s')$. The VFM property requires that $J^*(s) \geq J^*(s')$. In words, all else being equal, it is advantageous to start at a state with a greater number of items available.
\end{property}

Before we describe the second property we find it convenient to introduce some notation. For $\alpha \in \Z_+^M$, define a mapping $\tilde{S}_{\alpha}: \Sscr
\rightarrow \Sscr$ according to $\tilde{S}_{\alpha}(s) = s'$ with
$t(s') = t(s)$ and $x(s')_m = (x(s)_m - \alpha_m)^+$ for all $m$.
$\tilde{S}_{\alpha}(s)$ is thus the state obtained if one were
permitted to employ some set of activities (which presumably
resulted in $X_{t} = \alpha$) but without incurring the use of a
time-step.

\begin{property}
\label{prop:ir}
\noindent\textbf{Immediate Rewards: } For all $\bar{s} \in \Sscr$ and $\alpha \in \times_m \{0,1,\dots,x(\bar{s})_m\}$,
$$J^*(\bar{s}) \leq
g\left(x(\bar{s}), x(\bar{s}) - \alpha,t(\bar{s})\right)
+
J^*\left(\tilde{S}_{\alpha}(\bar{s})\right).$$
\end{property}

This property states that it is advantageous if one were able to deplete items
\emph{without incurring the use of a time-step}. In particular, if instead of
starting at some state $\bar{s} \in \Sscr$, one started
at state $\tilde{S}_{\alpha}(\bar{s})$ and was in addition given reward for the
depletion of $\alpha$ items, this property requires that the value of the second scenario be at least as large as the first.

The two properties we have developed thus far for the optimal value
function $J^*$ are what we fundamentally need to prove an approximation guarantee for the myopic heuristic. We now present the proof of our main approximation guarantee which assumes the VFM and IR properties.

We first introduce relevant notation: For $\alpha \in \mathbb{Z}^M_+$, define a mapping
$\hat{S}_{\alpha}: \Sscr \rightarrow \Sscr$ that specifies the next
state obtained if one employed an activity set in state $s$ which
resulted in  $X_{t,m} = \alpha_m \wedge x(s)_m$ for all $m$. In
particular, we define $\hat{S}_\alpha$ according to
$\hat{S}_{\alpha}(s) = s'$ with $t(s') = t(s)+1$ and $x(s')_m =
(x(s)_m - \alpha_m)^+$ for all $m$.


\begin{theorem}
\label{theorem:bound1} Assuming $J^*$ satisfies Properties \ref{prop:vfm} and \ref{prop:ir}, we have for all $s \in \Sscr$,
$\frac{J^*(s)}{J^{\pi^g}(s)} \leq 2.$
\end{theorem}

\begin{proof}
The proof proceeds by induction on the number of time steps that
remain in the horizon, $T-t(s)$. The claim is trivially true if
$t(s) = T-1$ since both the myopic and optimal policies coincide in
this case. Consider a state $s$ with $t(s) < T-1$ and assume the
claim true for all states $s'$ with $t(s') > t(s)$.

Now if $\pi^*(s) = \pi^{g}(s)$ then the next states encountered in both systems are
identically distributed so that the induction hypothesis immediately
yields the result for state $s$. Consider the case where $\pi^*(s) \neq \pi^g(s)$. Denote by $X^*_{t(s)}$ and $X^g_{t(s)}$ random vectors of depleted items in period $t(s)$ under optimal and myopic policies respectively at state $s$. Let $\mathbf{0}$ be an $M$ dimensional $0$ vector. We have:

\begin{equation}
\label{eq:opt}
\begin{split}
J^*(s|X^*_{t(s)}, X^g_{t(s)}) & = E[R(s,\pi^*(s))|X^*_{t(s)}] +
J^*(\hat{S}_{X^*_{t(s)}}(s))
\\
&
\leq E[R(s,\pi^*(s))|X^*_{t(s)}] + J^*(\hat{S}_{\mathbf{0}}(s))
\\
&
\leq
E[R(s,\pi^*(s))|X^*_{t(s)}] +
g(x(s), x(s) - X^g_{t(s)},t(s)+1)
+
J^*(\tilde{S}_{X^g_{t(s)}}(\hat{S}_{\mathbf{0}}(s)))
\\
&
\leq
E[R(s,\pi^*(s))|X^*_{t(s)}] +
g(x(s), x(s) - X^g_{t(s)},t(s))
+
J^*(\tilde{S}_{X^g_{t(s)}}(\hat{S}_{\mathbf{0}}(s)))
\\
&
=
E[R(s,\pi^*(s))|X^*_{t(s)}] + E[R(s,\pi^g(s))|X^g_{t(s)}]
+
J^*(\tilde{S}_{X^g_{t(s)}}(\hat{S}_{\mathbf{0}}(s)))
\\
&
=
E[R(s,\pi^*(s))|X^*_{t(s)}] + E[R(s,\pi^g(s))|X^g_{t(s)}]
+
J^*(\hat{S}_{X^g_{t(s)}}(s))
\\
&
\leq
E[R(s,\pi^*(s))|X^*_{t(s)}] + E[R(s,\pi^g(s))|X^g_{t(s)}]
+
2 J^{\pi^g}(\hat{S}_{X^g_{t(s)}}(s))
\end{split}
\end{equation}

where the first inequality follows from the assumed VFM property for $J^*$ upon noting that $x(\hat{S}_{X^*_{t(s)}}(s)) \leq x(\hat{S}_{\mathbf{0}}(s))$.
The second inequality follows from the IR property assumed of $J^*$ upon taking $\bar{s} = \hat{S}_{\mathbf{0}}(s)$ and $\alpha = X^g_{t(s)}$. The third inequality follows from Assumption \ref{as:g} since $g$ was assumed non-increasing in time. The third equality follows from the identity
$\tilde{S}_{X^g_{t(s)}}(\hat{S}_{\mathbf{0}}(s)) =  \hat{S}_{X^g_{t(s)}}(s)$ which in turn is simply a consequence of the definitions of $\tilde{S}_\alpha$ and $\hat{S}_\alpha$.
The final inequality follows from the induction hypothesis.

Now,
$$
J^{\pi^g}(s)
=
E\left[R(s,\pi^g(s)) + J^{\pi^g}(\hat{S}_{X^g_{t(s)}}(s))\right]
$$
and $E[R(s,\pi^g(s))] \geq E[R(s,\pi^*(s))]$ by the definition of the myopic policy $\pi^g$ so that taking expectations in \eqref{eq:opt}, we have:
\[
\begin{split}
J^*(s)
&
=
E[J^*(s|X^*_{t(s)}, X^g_{t(s)}) ]
\\
&
\leq
E[R(s,\pi^*(s))] + E[R(s,\pi^g(s))]
+
2 E[J^{\pi^g}(\hat{S}_{X^g_{t(s)}}(s))]
\\
&
\leq
2 J^{\pi^g}(s)
\end{split}
\]
This concludes the proof.
\qed
\end{proof}


\subsubsection{Performance with an approximate Myopic Oracle: } We will subsequently encounter a number of examples for which the set $\Ascr$ is exponentially large, but admits some implicit polynomial representation allowing for efficient solutions to the myopic problem
$$
\max_{A \in \Ascr} E[R(s,A)].
$$
Sometimes, however, this problem may itself be difficult to solve. In such scenarios the use of an oracle that is an $\alpha$-approximation to this subproblem is in fact a $1+\alpha$-approximation to the original stochastic depletion problem. In particular, assume $\pi^{\rm approx}:\Sscr \rightarrow \Ascr$ satisfies
$$
E[R(s,\pi^{\rm approx}(s))] \geq \frac{1}{\alpha}\max_{A \in \Ascr} E[R(s,A)].
$$
for all $s \in \Sscr$. One may then establish the following result whose proof  is omitted but entirely analogous to Theorem \ref{theorem:bound1} above:
\begin{theorem}
\label{theorem:bound2}
Assuming $J^*$ satisfies Properties \ref{prop:vfm} and \ref{prop:ir}, we have for all $s \in \Sscr$,
$\frac{J^*(s)}{J^{\pi^{\rm approx}}(s)} \leq 1 + \alpha.$
\end{theorem}

\section{Families Satisfying the VFM and IR Properties}
\label{sec:families}

The previous section identified two abstract properties - namely, the VFM and IR properties, that if satisfied yield uniform performance loss guarantees for the myopic heuristic, via Theorems \ref{theorem:bound1} and \ref{theorem:bound2}. These properties are in general difficult to check. We establish in this section two simple yet fairly general families of stochastic depletion problems that satisfy Properties \ref{prop:vfm} and \ref{prop:ir} thereby guaranteeing that the myopic heuristic is a 2-approximation algorithm for those families. Although, there may certainly be other families of problems satisfying the VFM and IR properties, the families we identify in this section accommodate a number of interesting applications which will be the focus of Sections \ref{sec:control} and \ref{sec:matroid}.

\subsection{Submodular Stochastic Depletion Problems}

We consider problems for which $\{P_{t}(A)\}$ is a $[0,1]^M$-valued stochastic process for all $A \in \Ascr$. Assuming one chooses action $A$ at time $t$, $X_{t,m}$ is a Binomial-$(x_{t,m},P_{t,m}(A))$ random variable that given $x_{t,m}$ and $P_{t,m}(A)$, is independent of the past and $X_{t,m'}$ for $m' \neq m$.
We assume submodular rewards. In particular, we assume $g(x_t,x_{t+1},t) = w(\bar{x} - x_{t+1}) - w(\bar{x} - x_{t})$, where $w:
\Z_+^M \rightarrow \R$ satisfies:
\begin{assumption}
\label{as:submod}
$w:\Z_+^M \rightarrow \R$ satisfies:
\begin{enumerate}
\item (Monotonicity) $w(y) \geq w(y')$ for $y \geq y'$.
\item (Submodularity) For $e \in \Z_+^M$, $w(y+e) - w(y) \leq w(y'+e) - w(y')$ if $y \geq y'$.
\end{enumerate}
\end{assumption}

Such a class of functions clearly satisfies Assumption \ref{as:g}. We need to demonstrate the VFM and IR properties. Recall that we will consider a clairvoyant optimal algorithm that knows a-priori the realizations of the sample paths of the $P_t$ processes. We first demonstrate the IR property. It turns out that doing so requires only the monotonicity of $w$; the sub-modularity of $w$ is not required for this property to hold.

\begin{lemma}(Immediate Rewards)\label{lemma:ub_opt} 
We have for submodular stochastic depletion problems, for all $\bar{s} \in \Sscr$ and $\alpha \in \times_m \{0,1,\dots,x(\bar{s})_m\}$,
$$J^*(\bar{s}) \leq
w\left(\bar{x} - x(\bar{s}) + \alpha\right)
-
w\left(\bar{x} -  x(\bar{s})\right)
+
J^*\left(\tilde{S}_{\alpha}(\bar{s})\right)$$.
\end{lemma}
\begin{proof}
Consider using the optimal policy starting at state $\bar{s}$, and let $S_T^*$ be the random state under this  policy at the end of the time horizon (that is, at time $T$), so that:
\begin{equation}
\label{eq:immediate_reward}
J^*(\bar{s})
=
E\left[w(\bar{x} - x(S_T^*))\right]
-
w(\bar{x} - x(\bar{s}))
\end{equation}
where the expectation is over the randomness in the system - namely,
the random item depletion defined by the $P_t$ sequences and chosen
activities. Similarly, let $\tilde{S}_T^*$ be the random state under the optimal policy at the end of the time horizon upon starting in state $\tilde{S}_{\alpha}(\bar{s})$ and as above, we note:
\[
J^*\left(\tilde{S}_{\alpha}(\bar{s})\right)
=
E\left[w(\bar{x} - x(\tilde{S}_T^*))\right]
-
w(\bar{x} - x(\bar{s})+\alpha)
\]

Let us re-consider the optimal policy starting at state $\bar{s}$
and in particular, let us partition the initial set of items into a
set of \emph{fictitious} and \emph{real} items; we assume that we begin with
$\alpha_m$ fictitious items of type $m$ and $x(\bar{s})_m - \alpha_m$
real items of type $m$. This partitioning serves purely as a labeling
of items and does not impact the evolution of the system in any
fashion. In particular, if at some point in time $t$, we have
$x^f_{t,m}$ and $x^r_{t,m}$ fictitious and real items of type $m$
respectively, then using activity set $A$ results in the depletion
of $X^f_{t,m}$ and $X^r_{t,m}$ fictitious and real items respectively
where $X^f_{t,m}$ is a Binomial-$(x^f_{t,m},P_{t,m}(A))$ random
variable and  $X^r_{t,m}$ is a Binomial-$(x^r_{t,m},P_{t,m}(A))$
random variable (so that $X^f_{t,m} + X^r_{t,m} \eqd X_{t,m}$), and
we are left with $x^f_{t,m} - X^f_{t,m}$ and $x^r_{t,m} - X^r_{t,m}$
fictitious and real items respectively.  Let $\sum_t X^f_{t,m}$
($\sum_t X^r_{t,m}$) denote the number of fictitious (real) items of
type $m$ depleted at the end of the time horizon by the optimal
policy starting in state $\bar{s}$.

We now make two critical observations:
\begin{enumerate}

\item We observe that $\sum_t X^f_{t,m} \leq \alpha_m$ for all $m$ by construction.

\item Due to the fact that given $s$ and a choice of activity, the depletion of a given item
of type $m$ at time $t$ is independent of the past and  the depletion of any
other item in the system at that time, $\sum_t X^r_{t,m}$ may be viewed as the
number of items of type $m$ depleted under some induced randomized
sub-optimal policy, say $\pi'$, starting at state
$\tilde{S}_{\alpha}(\bar{s})$.
This induced policy $\pi'$ assumes in state $\tilde{S}_{\alpha}(\bar{s})$ the existence of an additional $\alpha_m$ items of each type $m$ and simulates depletion of those items without garnering any reward for them. It operates like the optimal policy would but on this modified state.
More specifically, letting $s_0 = \tilde{S}_{\alpha}(\bar{s})$ and $\alpha_0 = \alpha$, we have:
\[
\pi'(s_0) = \pi^*((x(s_0) + \alpha_0,t(s_0)))
\]
Defining $\alpha_1 = x(S((\alpha_0,t(s_0)),\pi'(s_0)))$ and $s_1 = S(s_0,\pi'(s_0))$,
\[
\pi'(s_1) = \pi^*((x(s_1) + \alpha_1,t(s_1)))
\]
In general, defining $\alpha_t = x(S((\alpha_{t-1},t(s_{t-1})),\pi'(s_{t-1})))$ and $s_t = S(s_{t-1},\pi'(s_{t-1}))$,
\[
\pi'(s_t) = \pi^*((x(s_t) + \alpha_t,t(s_t)))
\]
It is worth noting that $\alpha_0 - \alpha_t = \sum_{t'=t(\bar{s})}^t X^f_{t'}$ while $x(s_0) - x(s_t) = \sum_{t'=t(\bar{s})}^tX^r_{t'}$.
\end{enumerate}

We consequently have:
\[
\begin{split}
J^*\left(\tilde{S}_{\alpha}(\bar{s})\right)
&
\geq
J^{\pi'}\left(\tilde{S}_{\alpha}(\bar{s})\right)
\\
&
=
E\left[w\left(\bar{x} -  x(\bar{s}) + \alpha + \sum_t X^r_{t}\right)\right]
-
w\left(\bar{x} - x(\bar{s})+\alpha\right)
\\
&
\geq
E\left[w\left(\bar{x} -  x(\bar{s}) + \sum_t X^f_t + \sum_t X^r_{t}\right)\right]
-
w\left(\bar{x} - x(\bar{s})+\alpha\right)
\\
&
=
E\left[w(\bar{x} - x(S_T^*))\right]
-
w\left(\bar{x} - x(\bar{s})+\alpha\right)
\\
&
=
J^*(\bar{s})
+
w(\bar{x} - x(\bar{s}))
-
w\left(\bar{x} - x(\bar{s})+\alpha\right)
\end{split}
\]
where the first inequality follows from the optimality of $\pi^*$
among all non-anticipatory policies. The first equality follows
from our definition of the policy $\pi'$ in Observation 2 and from
the definition of $\tilde{S}_{\alpha}(\bar{s})$. The second
inequality follows from the monotonicity of the function $w$ and
Observation 1: $\sum_t X_t^f \leq \alpha$. The second equality is
again by our construction of the $X^f_t$ and $X^r_t$ processes. The
final equality follows from \eqref{eq:immediate_reward}. This
completes the proof. \qed
\end{proof}

While the IR property required only the monotonicity of $w$, the VFM property requires both the monotonicity of $w$ as also its submodularity. This result is intuitive: a controller that starts at state $s$ may simply assume
that it starts at state $s'$ and track state evolution accordingly.
Assuming submodular rewards, applying the optimal policy to this
(incorrectly tracked) state trajectory guarantees the policy a
total expected reward of at least $J^*(s')$, so the optimal
policy must certainly do at least as well. The submodularity required is somewhat subtle, but it is simple to construct counterexamples in the absence of submodularity. We have:

\begin{lemma}(Value Function Monotonicity)\label{lemma:mono_greedy}
We have for Submodular Stochastic Depletion problems, for all $s,s' \in \Sscr$ s.t. $x(s) \geq x(s'), t(s) = t(s')$,
$J^*(s) \geq J^*(s').$
\end{lemma}
\begin{proof}
Consider a coupling of the systems starting at state $s$ and $s'$
wherein both systems witness identical sample paths for the item
depletion processes defined by the $P_t$ sequences. More precisely, assuming that at
time $t$, the systems are in states $s_t$ and $s'_t$ respectively,
then given $P_t$, the number of items depleted in both systems are
coupled so that if $x(s_t) \geq x(s'_t)$, and we employ activity set
$A$ in both systems, then, for all $m$, the number of successfully
depleted items of type $m$ in the $s_t$ system, $X_{t,m}$ ($\sim
{\rm Binomial}(x_{t,m},P_{t,m}(A)))$  and the number of
successfully depleted items of type $m$ in the $s'_t$ system,
$X'_{t,m}$ ($\sim {\rm Binomial}(x'_{t,m},P_{t,m}(A)))$ satisfy
$X_{t,m} = X'_{t,m} + Y_{t,m}$ where $Y_{t,m}$ is an independent
{\rm Binomial}-$(x(s_t)_m - x(s'_t)_m, P_{t,m}(A))$ random variable. A
symmetric situation must hold if $x(s'_t) \geq x(s_t)$.

Now assume that the system starting at $s'$ uses an optimal policy whereas the system starting at state $s$ mimics the actions of the $s'$ system (call this policy $\bar{\pi}$). It is simple to see that $\bar{\pi}$ is an admissible non-anticipatory policy.

Under our coupling, we have at $t=t(s)$, that the number of items of type $m$ depleted in the system starting at state $s$ is greater than the number of items depleted in the system starting at state $s'$. That is,  $X_{t,m} \geq X'_{t,m}$. It then follows that,
\[
\begin{split}
R(s,\bar{\pi}(s))
&
=
w(\bar{x} - x(s) + X_t)
-
w(\bar{x} - x(s))
\\
&
\geq
w(\bar{x} - x(s) + X'_t)
-
w(\bar{x} - x(s))
\\
&
\geq
w(\bar{x} - x(s') + X'_t)
-
w(\bar{x} - x(s'))
\\
&
=
R(s',\pi^*(s')).
\end{split}
\]
That is, the reward earned in the system starting at state $s$ is
higher than that in the system starting at state $s'$; the first
inequality above uses the monotonicity of $w$, the second inequality
employs the submodularity of $w$. Now in addition, by our coupling, both systems transition to states $S_{t(s)+1}$ and $S'_{t(s)+1}$ respectively satisfying $x(S_{t(s)+1}) = x(s) - X_t = (x(s') - X'_t) + (x(s) - x(s') - Y_t) \geq x(s') - X'_t = x(S'_{t(s)+1})$, so that we may repeat the above argument for time $t(s)+1$.
Continuing in this fashion we see that in every time step, the
$\bar{\pi}$ controlled system starting at state $s$ earns at least
as large a reward as the $\pi^*$ controlled system starting in state
$s'$. Taking expectations over the random item depletions (i.e. the $X_t$ and $X'_t$ processes), we have
$J^{\bar{\pi}}(s) \geq J^*(s')$. Since $J^*(s) \geq
J^{\bar{\pi}}(s)$, we are done.
\qed
\end{proof}

In light of Lemmas \ref{lemma:ub_opt} and \ref{lemma:mono_greedy}, Theorem \ref{theorem:bound1} lets us conclude that the myopic heuristic is a $2$-approximation algorithm for Submodular Stochastic Depletion problems.

\subsection{Linear Decaying Stochastic Depletion Problems}

We consider here a family of stochastic depletion problems closely related to the family just considered that also admit the VFM and IR properties and are particularly useful for many applications. As before, we consider problems for which $\{P_{t}(A)\}$ is a $[0,1]^M$-valued stochastic process for all $A \in \Ascr$. Assuming one chooses action $A$ at time $t$, $X_{t,m}$ is a Binomial-$(x_{t,m},P_{t,m}(A))$ random variable that given $x_{t,m}$ and $P_{t,m}(A)$, is independent of the past and $X_{t,m'}$ for $m' \neq m$.
We assume linear rewards that are non-increasing in time. In particular, we assume $g(x_t,x_{t+1},t) = \sum_m w_{m,t}(x_{t,m} - x_{t+1,m})$ where $w_{m,t}$ is a non-negative, non-increasing function of $t$ for all $m$ (for the special case where $w_{t,m} = w_m \geq 0$ for all $t$ this is merely a special case of a submodular stochastic depletion problem model we have considered). We can verify the immediate rewards property for such systems via a proof that closely follows Lemma \ref{lemma:ub_opt} and that may be found in the appendix:

\begin{lemma}(Immediate Rewards)\label{lemma:ub_opt_ld} 
We have for Linear Decaying Stochastic Depletion problems, for all $\bar{s} \in \Sscr$ and $\alpha \in \times_m \{0,1,\dots,x(\bar{s})_m\}$,
$$J^*(\bar{s}) \leq
\sum_m \alpha_mw_{t(\bar{s}),m}
+
J^*\left(\tilde{S}_{\alpha}(\bar{s})\right),$$
\end{lemma}

In addition, we may verify the the VFM property. The proof of the following Lemma is essentially identical to that of Lemma \ref{lemma:mono_greedy} and is omitted.

\begin{lemma}(Value Function Monotonicity)\label{lemma:mono_greedy_ld}
We have for Linear Decaying Stochastic Depletion problems, for all $s,s' \in \Sscr$ s.t. $x(s) \geq x(s'), t(s) = t(s')$,
$J^*(s) \geq J^*(s').$
\end{lemma}

In light of Lemmas \ref{lemma:ub_opt_ld} and \ref{lemma:mono_greedy_ld}, Theorem \ref{theorem:bound1} lets us conclude that the myopic heuristic is a $2$-approximation algorithm for Linear Decaying Stochastic Depletion problems.

\subsection{A Worst Case Example}
Having established the VFM and IR properties for the two families of stochastic depletion problems just discussed, we immediately have that the myopic policy has expected value within a factor of $2$ of the optimal policy for problems from either family. This analysis is sharp. In particular, we now present a problem instance that is in fact a member of both problem families and for which we have that the optimal policy has expected value that is a factor of $2-\eps$ better than the myopic policy; $\eps > 0$ can be made arbitrarily small.

\begin{example}(Myopic Sub-Optimality)
\label{ex:suboptimality}
Consider the case where $M = 2$ and $T = 2$, $g(x,x',t) = (x_1 - x'_1) + (1-\eps)(x_2 - x'_2)$. Assume that $\bar{x}_1 = \bar{x}_2 = 1$ and that $x_{0,1} =
x_{0,2} = 1$. Let $\Ascr =\{1,2\}$. The (deterministic) $P_t$ processes are defined as:
\[
\begin{split}
{\rm For} \ A=1: \
P_{0,1}(1) =  1, P_{0,2}(1) = 0, P_{1,1}(1) = 1, P_{1,2}(1) = 0.
\\
{\rm For} \ A=2: \
P_{0,1}(2) =  0, P_{0,2}(2) = 1, P_{1,1}(2) = 0, P_{1,2}(2) = 0.
\end{split}
\]

In words, the item of type $m=1$ may be depleted in either time step via the use of $A = 1$, whereas the item of type $m=2$ may be depleted only in the first time-step via $A = 2$. Only one of activity $1$ or $2$ may be employed within a given time-step.
The myopic heuristic will first choose activity set $\{1\}$ (which earns a reward of $1$ via the depletion of the $m=1$ type item) over activity set $\{2\}$  (which earns a reward of $1-\eps$ via the depletion of the $m=2$ type job). Consequently, under the myopic heuristic, $x_{1,1} = 0, x_{1,2} = 1$ and the heuristic is unable to complete the one remaining job in the second time step, earning a total reward of $1$.
An optimal schedule would first choose activity set $\{2\}$ (which earns a reward of $1-\eps$ via the completion of the $m=2$ type job). Consequently, under the optimal schedule, $x_{1,1} = 1, x_{1,2} = 0$ and the heuristic is able to complete the one remaining job in the second time step via the use of activity set $\{1\}$ earning a total reward of $2-\eps$. We thus see that $J^*(s_0) = (2-\eps) J^{\pi^g}(s_0)$ here.
\end{example}

\section{Applications: Stochastic Control}
\label{sec:control}
In the previous section, we presented two families of stochastic depletion problems for which the myopic heuristic is a $2$-approximation algorithm. We now consider several problems of stochastic control that are easily seen to be members of these families. We thereby establish uniform performance guarantees for myopic policies for these stochastic control problems.

\subsection{Service Policies for Simple Queueing Models} The following is a discrete time version of a `call-center' queueing model that has received a good deal of recent attention: We have $I$ buffers and $J$ servers. Each buffer sees a general discrete time arrival process with the restriction that a given buffer can see at most a single arrival in a given time slot. For example, each buffer $i$, may see an independent  ${\rm Bernoulli}(\lambda_i)$ arrival process. A given server $j$ may be used to service any single job in the system in a given time slot. In particular, should server $j$ be used to service a job arriving to buffer $i$, the service time is assumed to be an independent geometric random variable with mean $\mu_{i,j}$ (possibly $\infty)$. We allow for pre-emption in our service discipline. Consider the following natural objective: completion of a job that arrives at buffer $i$ earns a non-negative reward $r_{i,d}$ where $d$ is the time that job has remained in the system (that is, the delay experienced by that job). We assume $r_{i,d}$ is non-incerasing in $d$. At every point in time one must decide on a matching between servers and available jobs with a view to maximizing the expected reward earned over $T$ periods.

It is not difficult to see that the above problem is an example of a Linear Decaying Stochastic Depletion Problem. In particular, we define an item type for every tuple $(i,\tau)$ where $i = 1,2,\dots, I$ and $\tau=0,1,\dots,T-1$. Thus an item type $m$ is associated with an arrival buffer $i_m$ and an arrival time $\tau_m$. We can have at most a single item of a given type, i.e. $\bar{x}_m = 1$. The set of feasible activities $\Ascr$ is simply the set of all matching of servers to item types. Given a particular matching, the probability of depletion for a given item type (or job), is determined by the server matched to that job or $0$ if no server is matched to it. Of course, a job may not be depleted prior to arrival. In particular, we have for item type $m=(i_m,\tau_m)$
\[
P_{t,m}(A) = \mathbf{1}_{t \geq \tau_m}\left(\sum_j  \mathbf{1}_{(i,j) \in A} 1/\mu_{i,j} \right)
\]
We define our reward function $g$ according to $g(x_t,x_{t+1},t) = \sum_m w_{m,t}(x_{t,m} - x_{t+1,m})$ where we assume $w_{m,t} = r_{i_m,(t-\tau_m)^+}$. In particular, the reward generated in the $t$th time step is given by $\sum_{i=1}^I \sum_{t'=0}^{t} r_{i,t-t'}X_{i,t-t'}$ where $X_{i,t-t'} = 1$ if a job arriving to buffer $i$ at the start of time step $t'$ was completed at time $t$, and is $0$ otherwise. Since both the VFM and IR properties hold for this family of stochastic depletion problems, we have via Lemmas \ref{lemma:ub_opt_ld} and \ref{lemma:mono_greedy_ld}, and  Theorem \ref{theorem:bound1} that the myopic policy generates total expected rewards that are within a factor of $2$ of the optimal policy. In fact, we have shown that this performance guarantee holds relative to an optimal policy that has full knowledge of the entire job-arrival process!

A natural continuous time variant of the problem above has been the subject of much recent study and results on \emph{throughput} optimality in a certain heavy traffic regime are available (see, for instance, \cite{Bassamboo_QSTA}); the formulation we have discussed focuses on a different objective and complements that body of work. It is also interesting to note the analogy of our myopic policy with the so-called $c-\mu$ scheduling rules (see for instance, \cite{VanMieghem}) for scheduling jobs arriving to multiple buffers served by a single server with a view to minimizing total delay cost (every job incurs a buffer dependent, typically linear, delay cost). In our model, we maximize reward as opposed to minimizing cost; rewards decrease with delay but are necessarily non-negative.

\subsection{Stochastic Broadcast Scheduling}
\label{ssec:broadcast}

We consider a broadcast communication system where a single data item may be simultaneously transmitted to multiple users. In particular, we consider the following problem:  we have a set of $U$ users (indexed by $u$) and a finite set of data items or `pages', $P = \{1,\ldots,n\}$, indexed by $i$. In every time slot $\tau \in \{1,\dots,T-1\}$, any given user may generate a request for some page (or pages) he has not requested in the past. We assume that every request for a page is associated with a deadline $d\in\{1,\ldots,T-1\}$. Should a request for page $i$ by user $u$ be successfully satisfied prior to its deadline, the transmitter earns non-negative reward $r_i^u$.  We assume that the arrival process governing requests from users, as also the deadlines associated with those requests are exogenous stochastic processes and further assume a (known) bound on the number of requests that may arrive in any given time slot. In each time slot, a single page can be transmitted (although in what follows we could as well consider allowing up to $b$ pages). Due to the broadcast nature of the system, this transmitted page may be transmitted simultaneously to up to $k \geq 1$ users. The communication channel to users is stochastic, so that should a page be transmitted to a particular user $u$ in time $t$, that user receives the page with some channel dependent probability $P_t^u$ which is itself an exogenous stochastic process. In each time slot, one must decide which page to transmit and to which $k$ users in order to maximize the expected reward accrued over $T$ time slots.

Approximation algorithms for \emph{deterministic} broadcast scheduling (where transmissions are successful with probability $1$ so that $P_t^u = 1$ for all $t,u$) have received quite a bit of attention. The best known approximation algorithm is a $4$-approximation due to \cite{BarNoy02}. Without any constraints on the number of requests which can be satisfied by a single broadcast (i.e. $k = \infty$), the best known algorithm is a $4/3$-approximation algorithm due to \cite{Gandhi02}. The best known online algorithm for the same is a $2$-approximation due to \cite{Kima04}. Specializing to this deterministic case, our myopic online algorithm improves upon the offline results for `finite batching' in \cite{BarNoy02} albeit for uniform item sizes. Modeling stochasticity in communication channels to users is important since in real world systems, congestion and various physical phenomena cause significant uncertainty in the successful transmission of pages. Scheduling communications over stochastic channels is, of course, the focus of a substantial body of work in communications engineering. See, for instance, \cite{srikant_ton_05, tassiulas_infocom_97,tassiulas_globecom_02} for models closely related to the broadcast scheduling model we have presented. Most of that body of work is either simulation driven or else focuses on coarser performance metrics (such as throughput optimality).

The stochastic broadcast scheduling problem we have presented may be cast
as a Linear Decaying Stochastic Depletion Problems. Every request is associated with four parameters $(u,i,\tau,d)$ representing the user, page, time of request and request deadline respectively. We associate an item type with each such request. Thus an item type $m$ is identified by a request by user $u_m$ for page $i_m$, with an arrival time, $\tau_m$, and deadline, $d_m$. An activity $A \in \Ascr$ is simply an assignment of a given page to $k$ users.  Given a particular choice of activity $A$, the probability of depletion of a particular item type is simply given by the quality of the channel to the user corresponding to that type, $P_t^u$, provided that user $u$ is served under activity $A$; else it is $0$. Of course a request may not be satisfied prior to its arrival or following the expiration of its deadline, so that for item type $m = (u_m,i_m,\tau_m,d_m)$:
\[
P_{t,m}(A) = \mathbf{1}_{d_m > t \geq \tau_m} \mathbf{1}_{(i_m,u_m) \in A} P_t^{u_m}
\]
We define our reward function $g$ according to $g(x_t,x_{t+1},t) = \sum_m w_m(x_{t,m} - x_{t+1,m})$ where we assume $w_m = r_{i_m}^{u_m}$. This is a Linear Decaying Stochastic Depletion Problems and since both the VFM and IR properties hold for this family of stochastic depletion problems, we have via Theorem \ref{theorem:bound1} that the myopic policy generates expected value that is within a factor of $2$ of the optimal policy.


It is interesting to consider a special case of the stochastic broadcast scheduling problem we have presented: In particular, assuming that all requests are known at time $t=0$ and further that all these requests have deadline $T$, the myopic policy is in fact optimal if the channel to each user is `static', i.e., $P_t^u = C_u$, for all $u$. This is established via an interchange argument; the proof of the following Lemma may be found in the appendix.

\begin{lemma}
\label{lemma:greedy_opt}
If $P^u_t$ is constant ($= C_u$) for all
$u$, then, $J^{\pi^g}(s_0) = J^*(s_0)$ for all $s_0 \in \Sscr$.
\end{lemma}

Lemma \ref{lemma:greedy_opt} allows one to interpret the myopic
heuristic for the general stochastic broadcast scheduling problem as one that at every time-step $t$, makes the simplifying assumption that all channels are \emph{static} with success probabilities
given by $P^u_t$ and that no further arrivals will be observed. This is, in fact, a common engineering design principle for scheduling over dynamic channels. For instance,
\cite{tsibonis_twc_05, dua_icc_06, dua_tmc_07, huang_wcom_05, dua_globecom_06, Chou_tmm_06}
all derive optimal scheduling policies for problems similar to the broadcast scheduling problem here under the assumption of a static channel and other simplifying assumptions. The hope is that in conjunction with frequent channel state re-estimation (that is, frequent re-estimation of channel success probabilities), the use of scheduling schemes so derived may prove to be a very effective heuristic. In addition to being simple to implement and typically fast in practice, such an approach is robust to errors in specifying channel dynamics. Lemma \ref{lemma:greedy_opt} and Theorem \ref{theorem:bound1} thus lend theoretical support to this popular design principle. In particular, one may simply design a scheduling scheme assuming a static channel; one then employs this scheme in tandem with repeated channel re-estimation. Put another way, simply accounting for channel \emph{state} suffices to obtain levels of performance close to optimal.

\subsection{Dynamic Product Line Design}
\label{ssec:product_line}

Consider a firm that is capable of producing an array of related products that may potentially be sold to one or more customer segments, each distinguished by its willingness to pay for various product features. For a variety of reasons (manufacturing capacity and cost, marketing capabilities, etc.), the firm may be constrained in the number of different products it is capable of simultaneously offering for sale. Further, external competition may impose limitations on the prices the firm can post for a given product. Faced with these restrictions, the firm must decide on a product line to offer with a view to maximizing revenues. This is the essence of product line `design' problems that have been extensively considered in the operations research and marketing literature. For instance, the classic third degree price discrimination model of \cite{Pigou20} forms the basis of design principles that center on explicit market segmentation (see \cite{FWM72}). Alternatively, assuming a model of customer preference for various product attributes, one may consider optimizing the attributes of products offered for sale so as to maximize revenues; customers `self-select' product types that are of greatest appeal in this case. A number of product line design problems of this type have been considered in past literature; \cite{Moorthy84, Kohli90, VanRyzinMahajan99, Hopp05} are a few examples. A common thread to this work, however, is their consideration of static models. In reality, demand shocks and demand seasonality make the optimal product line design problem an inherently dynamic one.
For instance, consider the following example that illustrates the importance of accounting for seasonality in demand:
\begin{example}
A firm may offer at most one of two products (`outdated' or `new') for sale at any epoch (over two successive sales epochs) to two distinct  consumer segments - `bargain hunters' and `early adopters'. Bargain hunters will purchase only the outdated product for $1$ dollar in the first  period with probability $1$ and will make no purchase in period $2$. `Early adopters' will purchase only the new product with probability $1$ in either epoch for $1+\eps$ dollars. Assume we begin with an equal number of consumers in both segments. It is clear that a product line selection strategy that accounts for seasonality (by delaying the introduction of the new product to the second period) will earn about twice the revenues earned by a myopic strategy over two sales epochs.
\end{example}

Motivated by the consideration of issues such as the above, we consider the following dynamic product line design problem: a firm is capable of offering products from some set $\Pscr$ and must at any point in time offer a subset of products $A \subset \Pscr$ with $|A| \leq k$. The firm's products are purchased by $I$ consumer segments and we let $x_{t,i}$ denote the size of the $i$th segment in the $t$th sales epoch. Assuming that the product line offered at time $t$ is $A$, any segment $i$ consumer present in the market at that time $t$ will purchase a product in $A$ with known probability $P_{t,i}(A)$. Such a sale garners the seller revenue $p_i$ and the consumer is lost to the system. We assume that the firm has modeled the dynamics of the $P_t$ processes and wishes to maximize expected revenues over $T$ sales epochs.

It is worthwhile discussing some of the salient features of our model
\begin{enumerate}
\item We model the fact that a customer from a given segment may satisfy his requirement by the purchase of one of several product types, and that the probability of him purchasing a particular product type is influenced by the entire array of product types offered. Such a model permits customer self-selection but under the restriction that all substitutes are offered at the same price; this is precisely the type of model considered by \cite{VanRyzinMahajan99}. Alternatively, one may view the model as assuming that the seller has a means of directly segmenting customers (as in \cite{FWM72}) and allowing for segment specific prices.
\item We allow for general models of demand seasonality. In particular, we make no assumptions on the dynamics of the $P_t$ processes. Further, we explicitly model the impact of current sales on future demands (`market saturation').
\item We assume that the product line designer has available an estimate of market size within each consumer segment. Such an assumption is potentially valid in several industries; see \cite{john_deere} for a practical discussion of this issue.
\item Prices for each consumer segment are fixed. In reality this may arise, for example, due to the need to align with prices offered by competitors.
\end{enumerate}

It is simple to cast the above model as a Linear Decaying Stochastic Depletion Problems. In particular we associate with each customer segment $i$ an item type $m$. Our set of activities $\Ascr = \{A: A \in \Pscr, |A|\leq k\}$, and depletion probabilities for item type $m$ are specified according to $P_{t,m}(A) = P_{t,i_m}(A)$. We define our reward function $g$ according to $g(x_t,x_{t+1},t) = \sum_m w_{m,t}(x_{t,m} - x_{t+1,m})$ where we assume $w_{m,t} = p_{i_m}$. This is a Linear Decaying Stochastic Depletion Problems, so that Lemmas \ref{lemma:ub_opt_ld} and \ref{lemma:mono_greedy_ld} with Theorem \ref{theorem:bound1} immediately tell us that a myopic policy generates expected revenues within a factor of $2$ of the optimal policy. From a managerial perspective, this suggests a robust recipe for dealing with demand shocks and seasonality: at every opportunity for product line update, one simply solves a static product line design problem with suitably revised estimates of the relevant customer demand model and market sizes. As such, much of the existing methodology for product line design can be brought to bear on the problem without a significant cost to optimality.

\section{Applications: Stochastic variants of submodular maximization problems over Matroids}
\label{sec:matroid}

In this section, we turn our attention to the use of the stochastic depletion framework as a useful \emph{stochastic} analogue to submodular maximization problems over simple matroids such as the cardinality matroid and the partition matroid. A number of hard deterministic optimization problems can frequently be reduced to problems of this nature and in doing so, finding good approximation algorithms for these problems is reduced to the task of finding a good oracle for the myopic sub-problem (for a number of recently considered problems of this type, see for instance, \cite{Goundan}). Our hope is to produce good approximation algorithms for useful stochastic variants of such problems. As an illustration, we will later consider an important stochastic generalization of the AdWords Assignment problem considered by \cite{Fleischer} and \cite{Goundan}.

Given a set $E$, let $\Uscr = 2^E$. A \emph{cardinality matroid} is a subset of $\Uscr$ of the type $\Mscr = \{F \subset E: |F| \leq k\}$, where $k$ is an integer. A partition matroid is a subset of $\Uscr$ of the type $\Mscr = \{F \subset E: |F \bigcap E_i| \leq k_i \ \forall i\}$ where we assume $E = \uniond_{i=0}^n E_i$ and that integers $k_i$ for $i=0,1,\dots,n$ are given.

Consider optimization problems of the form
\begin{equation}
\label{eq:submod}
\max_{A \in \Mscr} f(A)
\end{equation}
where $f:2^E \rightarrow \R^+$ is a non-decreasing, submodular function. A number of interesting combinatorial optimization problems are reduced to such maximization problems where $\Mscr$ is a cardinality or partition matroid. We begin with establishing how such deterministic optimization problems are captured within the stochastic depletion framework.

\textbf{$\Mscr$ is a cardinality matroid: } We reduce \eqref{eq:submod} to a submodular stochastic depletion instance assuming $\Mscr$ is a cardinality matroid: We are given $M = |E|$ item types and assume that we begin with a single item of each type; i.e. $\bar{x}_m = x_{0,m} = 1$. Let $\Ascr = \{1, 2, \dots. |E|\}$ where $P_{t,m}(j) = 1$ if $m=j$ and $0$ otherwise. That is, we define an item type for each element of $E$ and in every time step we are allowed to deplete at most one item. We select as our reward function $g(x_t,x_{t+1},t) = f(\bar{x} - x_{t+1}) - f(\bar{x} - x_{t})$ and set $T=k$. Observe that the value of an optimal solution to this problem is precisely $J^*(\bar{x},0)$. With Lemmas \ref{lemma:ub_opt} and \ref{lemma:mono_greedy}, Theorem \ref{theorem:bound1} then immediately yields:

\begin{corollary}
The myopic heuristic is a $2$-approximation algorithm for maximizing a non-decreasing submodular function $f$ over a cardinality matroid.
\end{corollary}

We remark that this is a weaker result than the well known optimal approximation ratio of $\frac{e}{e-1}$ due to \cite{Nemhauser}. The analysis of Theorem \ref{theorem:bound1}, applies to a far broader class of problems, and in light of Example \ref{ex:suboptimality}, we can not expect a tighter guarantee for the greedy heuristic via that general line of analysis.

\textbf{$\Mscr$ is a partition matroid: } We reduce \eqref{eq:submod} to a submodular stochastic depletion instance assuming $\Mscr$ is a partition matroid: We are given $M = |E|$ item types and set $\Ascr = \{1, 2, \dots. |E|\}$. We index the elements of $E$ by $m$ which identifies a particular element of $E$ with a particular item type and assume that the first $|E_0|$ elements correspond to the elements of $E_0$, the next $|E_1|$ elements to the elements of $|E_1|$ and so forth. We set the time horizon $T = \sum_i k_i$ and define $n+1$ partitions of this horizon according to $T_j = [\sum_{l=0}^{j-1}k_l ,\sum_{l=0}^{j}k_l -1]$. We assume $P_{t,m}(j) = 1$ iff $m=j$ and $t \in \{T_j: m \in T_j\}$. We select as our reward function $g(x_t,x_{t+1},t) = f(\bar{x} - x_{t+1}) - f(\bar{x} - x_{t})$. In words, we define an item type for each element of $E$ and identify each subset $E_i$ with a partition of time. At any point in time $t \in T_i$, we are allowed to deplete at most one available item from the partition $E_i$. Observe that the value of an optimal solution to this problem is precisely $J^*(\bar{x},0)$. The myopic heuristic for this stochastic depletion problem corresponds precisely to the `local greedy' heuristic introduced by \cite{Fisher}, and we re-capture their result, namely:
\begin{corollary}
The local greedy heuristic is a $2$-approximation for maximizing a non-decreasing submodular function $f$ over a partition matroid.
\end{corollary}

Both classes of problems alluded to above have natural stochastic generalizations. As a simple example, one may consider a stochastic generalization to the problem of submodular maximization over a cardinality matroid which we refer to as the `stochastic selection problem'; as opposed to selecting at most $k$ elements from $E$, one is allowed $k$ attempts at  selecting elements of $E$. If at the $t$th selection attempt one attempts to select element $e \in E$, the attempt is successful with probability $P_t^e$ where $\{P^e_t\}$ is an arbitrary $[0,1]$ valued sequence specified for every $e \in E$. We would like to find an adaptive item selection policy that maximizes the expected value of the set of successfully selected items. It is easy to see that the stochastic selection problem includes as special cases appropriate stochastic generalizations of problems such as the maximum coverage problem.
The problem of adaptively selecting items so as to maximize the expected value of the set of successfully selected items is seen to be a submodular stochastic depletion problem using precisely the reduction for the cardinality matroid above and one immediately has the following result.

\begin{lemma}
The myopic heuristic is a $2$-approximation for the stochastic selection problem.
\end{lemma}

As an aside we note that if in addition, one assumes that $P_t^e = C$, a constant for all $t$ and $e$, it is simple to demonstrate that the myopic heuristic is, in fact, an $\frac{e}{e-1}$ approximation. This may be demonstrated as a corollary to the original result of \cite{Nemhauser}: one simply considers coupling the optimal and myopic schemes so that on each sample path, both schemes have an identical number of successful placements.

We now consider in some detail, a practically relevant stochastic generalization of the AdWords Assignment problem (\cite{Fleischer}). The deterministic problem may be reduced to the maximization of a submodular function over a partition matroid (see \cite{Goundan}). We reduce our stochastic generalization to a submodular stochastic depletion problem.

\subsection{Cost-per-Click AdWords Assignment}

Consider the following optimization problem faced by firms that serve ads on the internet. We are given a set of $N$ advertisers (indexed by $i$) and $K$ keywords (indexed by $k$). The $i$th advertiser has a budget $B_i (\geq 0)$ and submits to the firm a valuation $v_{i,k}$ for every keyword $k$. In every one of $T$ periods, a keyword from the set of $K$ keywords arrives according to some exogenous stochastic process. We assume that at most $C$ advertisers' ads can be assigned to the arriving keyword. We denote by $k_t$ the index of the keyword arriving at time $t$. Should an advertiser $i$ be assigned to an arriving keyword, $k_t$, at time $t$ and if in addition his ad is clicked on, he pays the firm  the minimum of $v_{i,k_t}$ and his remaining budget at time $t$; this payment is subtracted from his available budget. If the ad is not clicked on then no payments are made. We assume that should advertiser $i$ be assigned to keyword $k_t$ at time $t$, his ad is clicked on with probability $P_t^{i,k_t}$. Letting $V_{t,i}$ denote the random payment thus made by an advertiser in the $t$th period, we are interested in devising an adaptive ad-to-keyword assignment scheme that maximizes $E\left[\sum_{i=1}^N (B_i \wedge \sum_t  V_{t,i})\right]$, that is, the expected revenues earned by the firm.

The above problem was considered in a deterministic offline setting by \cite{Fleischer} where it was assumed that $P_t^{i,k} = 1$ for all $t,i,k$ and in addition the sequence $\{k_t\}$ of arriving keywords was specified a-priori; the variant we consider here is an important generalization to that model since in practice advertisers make payments only if their displayed ads are clicked on, which happens with some positive, but small, probability. In addition, our formulation also allows us to capture exogenous advertiser arrivals and departures from the system - in particular, we simply assume $P_t^{i,k} = 0$ for all times $t$ prior to a customers arrival to the system and following his departure from the system.

This problem is easily cast as a submodular stochastic depletion problem: In
particular, we define an item type $m$ for every advertiser-keyword-time triple $(i_m,k_m,t_m)$, and
assume a single item of each type, i.e. $x_{0,m} = \bar{x}_m = 1$ for all $m$. The set of feasible activities, $\Ascr$, is the set of all subsets of item types, such that each subset has cardinality at most $b$ and contains at most one type specific to a given user $i$. The probability that an item of type $m$ is depleted at time $t$ assuming one selects activity $A$ is given by,
$$
P_{t,m}(A) = \mathbf{1}_{\{(i_m,k_m,t_m) \in A, t = t_m, k_t = k_m\}}P_t^{i_m,k_m}.
$$
Finally, the reward function, $g(x_t,x_{t+1},t) = w(\bar{x}-x_{t+1}) - w(\bar{x} -x_t)$ where $w: \Z^M \rightarrow \R^+$ is defined according to
$$
w(x) = \sum_i  \left(B_i \wedge \sum_{m: i_m = i} v_{i_m,k_m} x_m\right)
$$
and thus satisfies Assumption \ref{as:submod}. We finally note
that the myopic sub-problem is trivial; it corresponds to choosing
the $b$ highest expected revenue advertisers. With Lemmas \ref{lemma:mono_greedy} and \ref{lemma:ub_opt}, Theorem
\ref{theorem:bound1} thus yields:
\begin{corollary}
The myopic heuristic is a $2$-approximation to Cost-per-Click AdWords Assignment
\end{corollary}

In our formulation, a feasible ads-to-keyword assignment was subject to a simple cardinality constraint: an arriving keyword could have at most $C$ ads assigned to it. We could instead consider using other, more complex constraints: in particular, in the formulation of \cite{Fleischer}, every ad is associated with a rectangle of a specific height and width, and every arriving keyword with an available rectangular display area; a feasible assignment of ads to keywords is
determined by a feasible packing of ad rectangles within the display rectangle. Using the max-weight rectangle packing $2+\eps$-approximation algorithm of \cite{Jansen} for the myopic subproblem yields via Theorem \ref{theorem:bound2}, a $3+\eps$-approximation guarantee for the myopic heuristic which matches the best known approximation guarantee available for the original deterministic problem (see \cite{Goundan}).

\section{Concluding Remarks}
\label{sec:conclusion}

We have in the present work introduced a general class of dynamic stochastic optimization problems - Stochastic Depletion problems. We believe this to be an interesting class of problems: in spite of being fairly general, stochastic depletion problems frequently admit a simple, near-optimal myopic control policy. This paper presented general conditions that guarantee the near optimality of a myopic control policy for a stochastic depletion problem and went on to verify these properties for broad families of stochastic depletion problems. This in turn yielded myopic approximation algorithms for a number of interesting dynamic optimization applications.

There are several directions that deserve continued study. From an algorithmic perspective, one may consider $k$-step lookahead policies as a generalization of the myopic ($1$-step lookahead) policies analyzed here. Such policies select, at every point in time, an action that is optimal for a problem with a horizon precisely $k$ time steps ahead. It would be interesting to understand whether, or under what conditions, such policies may be expected to dominate the myopic policy.

In addition to the applications in Section \ref{sec:control} and \ref{sec:matroid}, it would be interesting to explore other dynamic stochastic optimization problems that may be studied either within our framework or perhaps slight modifications to it. For instance, the generalized assignment problem (\cite{ShmoysGap}) is known to reduce to the maximization of a submodular function over a partition matroid. An interesting stochastic generalization of this problem that would allow for a number of interesting applications would involve making the successful placement of an item in a bin stochastic. Unfortunately, this particular stochastic generalization does not reduce to a stochastic depletion problem but is nonetheless very similar to one.

Another broad issue is identifying other families of stochastic depletion problems that satisfy the VFM and IR properties, or in another direction, identifying conditions under which we may not expect one of those properties to be satisfied. Yet another issue is the optimality of our approximation schemes: for deterministic variants of several of the application problems considered in this work such as the submodular maximization problems over matroids, there exist (typically, fairly complex) \emph{offline} algorithms that admit an approximation ratio of $\frac{e}{e-1}$ (see \cite{Calinescu}). This guarantee is known to be optimal. That is, no efficient approximation algorithm with a superior guarantee exists, unless $P=NP$. It would be interesting to understand whether an approximation ratio of $2$ is optimal in some sense for either of the families of problems for which we have established that guarantee in this paper. In any case, given that the best approximation guarantee we may expect is a factor of $\frac{e}{e-1}$, it is remarkable that a simple myopic scheme comes so close to achieving that guarantee and that this guarantee may be established in the generality of the stochastic depletion framework.

\section*{Acknowledgments}
The authors are grateful to Retsef Levi for his input over the course of this research - particularly the connections with matroid problems. Ben Van Roy pointed out potential connections between call-center scheduling and the extant combinatorial scheduling literature to the second author several years ago. The second author thanks Andreas Schulz for pointing out his interesting work on submodular maximization.

\small{
\bibliographystyle{ormsv080} 
\bibliography{MOR_1-16-08}
}

\begin{appendix}
\section{Miscellaneous Technical Proofs}
\begin{l3}(Immediate Rewards)
We have for Linear Decaying Stochastic Depletion problems, for all $\bar{s} \in \Sscr$ and $\alpha \in \times_m \{0,1,\dots,x(\bar{s})_m\}$,
$$J^*(\bar{s}) \leq
\sum_m \alpha_mw_{t(\bar{s}),m}
+
J^*\left(\tilde{S}_{\alpha}(\bar{s})\right),$$
\end{l3}
\begin{proof}
Consider using an optimal policy starting at state $\bar{s}$. Let us partition the initial set of jobs into a set of `fictitious' and `real' jobs; we assume that we begin with $\alpha_m$ fictitious jobs of type $m$ and $x(\bar{s})_m - \alpha_m$ real jobs of type $m$. This partitioning serves purely as a labeling of jobs and does not impact the system in any fashion. In particular, if at some point in time $t$, we have $x^f_{t,m}$ and $x^r_{t,m}$ fictitious and real jobs of type $m$ respectively, then using activity set $A$ results in the completion of $X^f_{t,m}$ and $X^r_{t,m}$ fictitious and real jobs respectively where $X^f_{t,m}$ is a Bernoulli-$(x^f_{t,m},P_{t,m}(A))$ random variable and  $X^r_{t,m}$ is a Bernoulli-$(x^r_{t,m},P_{t,m}(A))$ random variable (so that $X^f_{t,m} + X^r_{t,m} \eqd X_{t,m}$). The revenues earned are $\sum_m X^f_{t,m}w_{t,m}$ and $\sum_m X^r_{t,m}w_{t,m}$ and we are left with $x^f_{t,m} - X^f_{t,m}$ and $x^r_{t,m} - X^r_{t,m}$ fictitious and real jobs respectively.

Denote by $J^{f,*}(\bar{s})$ the expected reward-to-go under an optimal policy starting at state $\bar{s}$ earned from the completion of fictitious jobs. Likewise, we define  $J^{r,*}(\bar{s})$ as the expected reward-to-go under an optimal policy starting at state $\bar{s}$ earned from the completion of real jobs.
Now, by construction, $J^{f,*}(\bar{s}) + J^{r,*}(\bar{s}) = J^*(\bar{s})$. Since at best our scheduling policy can exhaust all fictitious jobs and since $w_{t,m}$ is non-increasing in $t$ for all $m$, $J^{f,*}(\bar{s}) \leq \sum_m \alpha_m w_{t(\bar{s}),m}$. Now, $J^{r,*}(\bar{s})$ may be viewed as the reward-to-go under some admissible policy $\pi$ starting at state $\tilde{S}_\alpha(\bar{s}))$. Noting that $x(\tilde{S}_\alpha(\bar{s}))_m$ is precisely the initial number of `real' jobs of type $m$, we then have:
$J^{r,*}(\bar{s}) = J^\pi(\tilde{S}_\alpha(\bar{s}))) \leq J^*(\tilde{S}_\alpha(\bar{s})))$. Consequently, we have,
\[
J^*(\bar{s}) = J^{f,*}(\bar{s}) + J^{r,*}(\bar{s}) \leq \sum_m \alpha_mw_{t(\bar{s}),m} + J^*(\tilde{S}_\alpha(\bar{s}))),
\]
which is the result.
\qed
\end{proof}

\begin{l5}
If $P^u_t$ is constant ($= C_u$) for all
$u$, then, $J^{\pi^g}(s_0) = J^*(s_0)$ for all $s_0 \in \Sscr$.
\end{l5}
\begin{proof}
Note that since the processes $P^u_t$ are deterministic here, we may without loss restrict attention to policies that are functions of only time $t$, and
$(x^1,x^2,\dots,x^U)$. Let $w^u_i = r_m$ for item $m = (i,u,\tau_m,d_m)$.  We define the set of myopic packets as
$\Pscr^*_0 = \argmax_j \sum_u C_uw^u_j x^u(s_0)_j$. Let us assume for
the sake of contradiction that in state $s_0$, no optimal policy transmits
a packet in $\Pscr^*_0$. Let $\pi^*$ be an optimal policy; $\pi^*(s_0) \notin
\Pscr^*_0$.

Define a policy $\bar{\pi}$ according to:
\begin{eqnarray*}
\bar{\pi}(s_0) =  \pi^g(s_0) &
\\
\bar{\pi}(s) =  \pi^*(f(s)) & \forall s \neq s_0.
\end{eqnarray*}
where $f:\Sscr \rightarrow \Sscr$ is defined according to $x^{\pi^g(s_0)}(f(s)) =  x^{\pi^g(s_0)}(s_0)$, $x^j(f(s)) = x^j(s) \ \forall j \neq \pi^g(s_0)$, and
$t(f(s)) = t(s)-1$.

Further, define a policy $\tilde{\pi}$ according to:
\begin{eqnarray*}
\tilde{\pi}(s) = & \pi^*(g(s)) & \forall s
\end{eqnarray*}
where $g:\Sscr \rightarrow \Sscr$ is defined according to $x(g(s))  = x(s), t(g(s)) =
t(s)+1$.

Let $\tau = \min \{t : \bar{\pi}(s_t) \in \Pscr^*_0\}$ (set $\tau$
to $\infty$ if the set is empty) and consider using policy
$\bar{\pi}$ for $t <\tau$ and policy $\tilde{\pi}$ thereafter. We
call this policy $\pi'$. Denote by $R^{\pi'}_t$ and
$R^{\pi^*}_t$ the random rewards earned in the $t$th time step under
the $\pi'$ and $\pi^*$ policies respectively, so that $J^*(s_0) =
E[\sum_t R^{\pi^*}_t]$ and $J^{\pi'}(s_0) = E[\sum_t R^{\pi'}_t]$.
Now observe that by our construction, $E[R^{\pi'}_0] >
E[R^{\pi^*}_0]$, $R^{\pi'}_t \eqd R^{\pi^*}_{t-1}$ for $0 < t
< \tau$ and $R^{\pi'}_t \eqd R^{\pi^*}_{t}$ for $t \geq \tau$. If
$\tau < \infty$, then it immediately follows that $J^*(s_0|\tau <
\infty) = J^{\pi'}(s_0|\tau<\infty)$; if $\tau = \infty$ then upon
noting that $E[R^{\pi'}_0] > E[R^{\pi^*}_{T-1}]$, we have
$J^*(s_0|\tau = \infty) > J^{\pi'}(s_0|\tau=\infty)$. Consequently,
$J^{\pi'}(s_0) \geq J^*(s_0)$, so that $\pi'$ is an optimal policy as well. This contradicts our assumption that no optimal policy transmits a packet in $\Pscr^*_0$ at $t=0$. We may thus assume without loss that an optimal policy transmits a packet in $\Pscr^*_0$ in the first time
step. This suffices for the proof.
\qed
\end{proof}

\end{appendix}

\end{document}